\documentclass[11pt,reqno]{amsart}

   \usepackage{amsthm}
   \usepackage{amsmath}
   \usepackage{amssymb}
   \usepackage{amsfonts}
   \usepackage{amscd}
   \usepackage[mathscr]{eucal}
   \usepackage{verbatim}
\allowdisplaybreaks

\newtheorem{theorem}{Theorem}[]
\newtheorem{lemma}[theorem]{Lemma}

\theoremstyle{remark}
\newtheorem{rem}[theorem]{Remark}

\newcommand{\C}{\mathbb C}
\newcommand{\Q}{\mathbb Q}

\newcommand{\ad}{\mathrm{ad}}
\newcommand{\Ad}{\mathrm{Ad}}
\newcommand{\GL}{\mathrm{GL}}

\newcommand{\GSp}{\mathrm{GSp}}
\newcommand{\PGSp}{\mathrm{PGSp}}
\newcommand{\SO}{\mathrm{SO}}
\newcommand{\SSp}{\mathrm{Sp}}
\newcommand{\gsp}{\mathfrak{gsp}}
\newcommand{\ssp}{\mathfrak{sp}}
\renewcommand{\sl}{\mathfrak{sl}}

\newcommand{\St}{\mathrm{St}}
\newcommand{\std}{\mathrm{std}}
\newcommand{\triv}{\mathbf{1}}

\newcommand{\OF}{\mathfrak{o}}
\newcommand{\p}{\mathfrak{p}}
\newcommand{\mat}[4]{{\setlength{\arraycolsep}{0.5mm}\left[
\begin{array}{cc}#1&#2\\#3&#4\end{array}\right]}}

\begin{document}

\title{On the adjoint $L$-function of the $p$-adic $\GSp(4)$}
\author[M. Asgari]{Mahdi Asgari}
\address{Department of Mathematics \\
Oklahoma State University \\
Stillwater, OK 74078-1058 \\ 
USA}
\email{asgari@math.okstate.edu}

\author[R. Schmidt]{Ralf Schmidt}
\address{Department of Mathematics  \\ 
University of Oklahoma \\
Norman, OK 73019-0315 \\
USA}
\email{rschmidt@math.ou.edu}

\begin{abstract}
We explicitly compute the adjoint $L$-function of those $L$-packets of representations 
of the group $\GSp(4)$ over a $p$-adic field of characteristic zero that contain non-supercuspidal 
representations. As an application we verify a conjecture of Gross-Prasad and Rallis in 
this case. The conjecture states that the adjoint $L$-function is holomorphic at $s=1$ if and 
only if the $L$-packet contains a generic representation.   
\end{abstract}

\maketitle

\section{Introduction}\label{intro}

Let $F$ be a non-archimedean local field of characteristic zero and let
$W'_F$ be the Weil-Deligne group of $F$. The conjectural local Langlands correspondence
for the group $\GSp(4,F)$ assigns to each irreducible admissible representation $\Pi$
of $\GSp(4,F)$ an $L$-parameter, i.e., an equivalence class of admissible
representations
$$
 \varphi_\Pi:\:W'_F\longrightarrow\GSp(4,\C).
$$
It was shown in \cite[Sect.\ 2.4]{RS} that there is a unique way to assign
$L$-parameters to the \emph{non-supercuspidal} irreducible, admissible
representations of $\GSp(4,F)$ such that certain desired properties of the
local Langlands correspondence hold. In this sense the local Langlands
correspondence is known for the non-supercuspidal representations of $\GSp(4,F)$;
see Table \ref{maintable} for a complete list of these representations.
In a few cases the $L$-packet of a non-supercuspidal representation is expected to also 
contain a supercuspidal representation.

The degree 4 and degree 5 $L$-factors resulting from the non-supercuspidal local Langlands
correspondence have been computed and tabulated in \cite[Tables A.8 and A.10]{RS}. 
In this article we treat the next smallest irreducible representation of the dual group,
namely the $10$-dimensional adjoint representation ${\rm Ad}$ of $\GSp(4,\C)$
on the complex Lie algebra $\ssp(4)$. Thus, given a non-supercuspidal, irreducible,
admissible representation $\Pi$ of $\GSp(4,F)$ with $L$-parameter $\varphi_\Pi$, we
compute
$$
 L(s,\Pi,{\rm Ad}):=L(s,{\rm Ad}\circ\varphi_\Pi).
$$
This is an easy calculation in most cases, but requires some arguments in a few.
The results are tabulated in Table \ref{polestable} below.

Having explicit formulas for all the adjoint $L$-functions, we immediately obtain 
the following case of a general conjecture of Gross-Prasad \cite[Conj.\ 2.6]{GP} 
and Rallis \cite[Prop.\ 5.2.2]{K} as a corollary; see Theorem \ref{GP-R} below\footnote{
The same statement is made in a preprint of W.~Gan and S.~Takeda 
\cite{gan-takeda} which became available after the completion of this work. It is 
very likely that their L-parameters coincide with ours.}.

\emph{Let $\Pi$ be a non-supercuspidal irreducible admissible 
representation of $\GSp(4,F)$. Then the $L$-packet of $\Pi$ contains a 
generic representation if and only if $L(s,\Pi,{\rm Ad})$ is holomorphic at $s=1$.} 

The analogous statement for $\GL(n,F)$ ``has been observed by many people'', loc cit. 
For a proof see \cite[Prop.\ 7.1]{JS}.

We note that there is some overlap between Theorem \ref{GP-R} and a result 
of  Jiang and Soudry.
In \cite{JS-annals} and \cite{JS} they attach to each admissible $L$-parameter
an irreducible, admissible representation of $\SO(2n+1,F)$ and prove that this
representation is generic if and only if its associated adjoint $L$-function
is holomorphic at $s=1$ \cite[Theorem 7.1]{JS}. In the special case $n=2$,
since $\SO(5,F)\cong\PGSp(4,F)$,
the representation of $\SO(5,F)$ corresponds to a representation of
$\GSp(4,F)$ with trivial central character. However, it is not immediately
clear that this version of the local Langlands correspondence coincides
with ours. To mention one difference, the Jiang-Soudry correspondence misses
those representations of $\GSp(4,F)$ whose central character is not a square,
since such representations are not a twist of a representation with trivial
central character. Also, the Jiang-Soudry correspondence does not assign
an $L$-parameter to the non-generic representations of type VIb
and XIb (see Table 1), both of which share an $L$-packet with a generic representation.

The authors would like to thank D.~Jiang and D.~Prasad for some helpful discussions. 
We are very grateful to J.~Cogdell for helpful discussions and feedback on an earlier 
version of this article. 

\section{Notation and Definitions}
\subsection{Group-theoretic definitions}
We realize the algebraic $\Q$-group $\GSp(4)$ as
$$
 \GSp(4)=\{g\in\GL(4):\:^tgJg=\lambda(g)J\text{ for some }\lambda(g)\in\GL(1)\},
$$
where
$$
 J=\begin{bmatrix}&&&1\\&&1\\&-1\\-1\end{bmatrix}.
$$
The kernel of the \emph{multiplier homomorphism} $g\mapsto\lambda(g)$ is by definition
the symplectic group $\SSp(4)$. The Lie algebra of $\SSp(4)$ is $10$-dimensional
and is given by
$$
 \ssp(4)=\{X\in\mathfrak{gl}(4):\:^tXJ+JX=0\}.
$$
Over the complex numbers, the Lie algebra of $\GSp(4)$ is a direct sum
$$
 \gsp(4)=\ssp(4)\oplus\mathfrak{z},\qquad
 \mathfrak{z}=\C\begin{bmatrix}1\\&1\\&&1\\&&&1\end{bmatrix}.
$$
The adjoint representation of $\GSp(4,\C)$ on $\gsp(4)$ preserves both summands  and, as 
representations of $\GSp(4,\C)$, we have 
\begin{equation}\label{liegspsspeq}
 \Ad_{\gsp}=\Ad_{\ssp}\oplus\triv. 
\end{equation}
We use $\Ad$ for $\Ad_{\ssp}$ in this article. 

The character lattice of $\SSp(4)$ is spanned by
\begin{equation}\label{characterbasiseq}
 e_1:\;\begin{bmatrix}a\\&b\\&&b^{-1}\\&&&a^{-1}\end{bmatrix}
 \longmapsto a\qquad\text{and}\qquad
 e_2:\;\begin{bmatrix}a\\&b\\&&b^{-1}\\&&&a^{-1}\end{bmatrix}
 \longmapsto b.
\end{equation}
We shall use the following generators for the root spaces in $\ssp(4)$.
\begin{alignat}{2}
 \label{rootspaces1eq}L_{e_1-e_2}&=\begin{bmatrix}0&1\\&0\\&&0&-1\\&&&0\end{bmatrix}
 &\qquad L_{-e_1+e_2}&=\begin{bmatrix}0\\1&0\\&&0\\&&-1&0\end{bmatrix}\\
 \label{rootspaces2eq}L_{e_1+e_2}&=\begin{bmatrix}0&&1\\&0&&1\\&&0\\&&&0\end{bmatrix}
 &\qquad L_{-e_1-e_2}&=\begin{bmatrix}0\\&0\\1&&0\\&1&&0\end{bmatrix}\\
 \label{rootspaces3eq}L_{2e_1}&=\begin{bmatrix}0&&&1\\&0\\&&0\\&&&0\end{bmatrix}
 &\qquad L_{-2e_1}&=\begin{bmatrix}0\\&0\\&&0\\1&&&0\end{bmatrix}\\
 \label{rootspaces4eq}L_{2e_2}&=\begin{bmatrix}0\\&0&1\\&&0\\&&&0\end{bmatrix}
 &\qquad L_{-2e_2}&=\begin{bmatrix}0\\&0\\&1&0\\&&&0\end{bmatrix}
\end{alignat}
The root system of $\SSp(4)$ is of type $C_2$,

\vspace{2ex}\qquad
\setlength{\unitlength}{4ex}
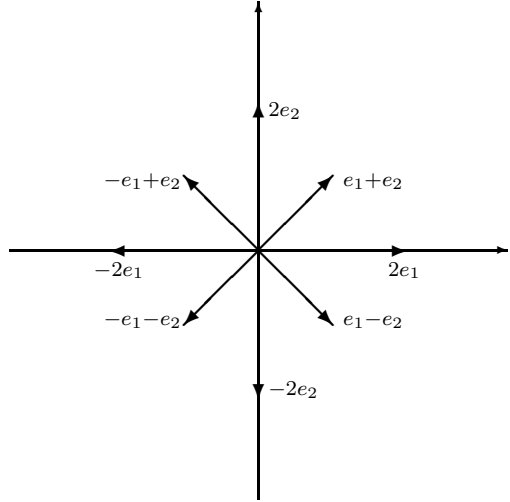
\begin{picture}(10,10)
 \put(0,5){\vector(1,0){10}}
 \put(5,0){\vector(0,1){10}}
 \thicklines
 \put(5,5){\vector(1,0){3}}
 \put(5,5){\vector(-1,0){3}}
 \put(5,5){\vector(0,1){3}}
 \put(5,5){\vector(0,-1){3}}
 \put(5,5){\vector(1,1){1.5}}
 \put(5,5){\vector(-1,-1){1.5}}
 \put(5,5){\vector(-1,1){1.5}}
 \put(5,5){\vector(1,-1){1.5}}
 \put(7.6,4.5){$\scriptstyle 2e_1$}
 \put(1.7,4.5){$\scriptstyle -2e_1$}
 \put(5.2,7.7){$\scriptstyle 2e_2$}
 \put(5.2,2.1){$\scriptstyle -2e_2$}
 \put(6.7,6.3){$\scriptstyle e_1+e_2$}
 \put(1.9,3.5){$\scriptstyle -e_1-e_2$}
 \put(6.7,3.5){$\scriptstyle e_1-e_2$}
 \put(1.9,6.3){$\scriptstyle -e_1+e_2$}
\end{picture}

The conjugacy classes of proper parabolic subgroups of $\GSp(4)$ are
represented by the minimal parabolic subgroup $B$, the Siegel parabolic
subgroup $P$, and the Klingen parabolic subgroup $Q$, consisting 
of matrices in $\GSp(4)$ of the following form, respectively :
$$
 B=\begin{bmatrix} *&*&*&*\\&*&*&*\\&&*&*\\&&&*\end{bmatrix},\qquad
 P=\begin{bmatrix} *&*&*&*\\ *&*&*&*\\&&*&*\\&&*&*\end{bmatrix},\qquad
 Q=\begin{bmatrix} *&*&*&*\\&*&*&*\\&*&*&*\\&&&*\end{bmatrix}.
$$
Setting 
\begin{equation}\label{Apdefeq}
 A'=\mat{}{1}{1}{}\,^t\!A^{-1}\mat{}{1}{1}{}\qquad\text{for }A\in\GL(2),
\end{equation} 
a typical element of $P$ can be written as $\mat{A}{*}{}{cA'}$ with
$c\in\GL(1)$ and $A\in\GL(2)$.
%%%%%%%%%%%%%%%%
\subsection{$p$-adic definitions}
Let $F$ be a non-archimedean local field of characteristic zero.  Let $\OF$ be its ring of
integers and $\p$ the maximal ideal of $\OF$.  We fix a generator
$\varpi$ of $\p$ once and for all.  A character $\chi$ of $F^\times$ is a
continuous homomorphism $F^\times\rightarrow\C^\times$.  It is unramified
if $\chi(\OF^\times)=\{1\}$.  A distinguished unramified character
is $\nu$, the normalized absolute value of $F$.  It has the property that
$\nu(\varpi)=q^{-1}$, where $q$ is the number of elements of the residue
class field $\OF/\p$.

We shall use the notation of \cite{ST} for representations of $\GSp(4,F)$
parabolically induced from one of the parabolic subgroups $B$, $P$ or $Q$.
If $\chi_1$, $\chi_2$ and $\sigma$ are characters of $F^\times$, then
$\chi_1\times\chi_2\rtimes\sigma$ denotes the representation of $\GSp(4,F)$
obtained via (normalized) parabolic induction from the character
$$
 \begin{bmatrix}a&*&*&*\\&b&*&*\\&&cb^{-1}&*\\&&&ca^{-1}\end{bmatrix}
 \longmapsto\chi_1(a)\chi_2(b)\sigma(c)
$$
of $B(F)$. If $\sigma$ is a character of $F^\times$ and $\pi$ is an
admissible representation of $\GL(2,F)$, we denote by $\pi\rtimes\sigma$ the
representation of $\GSp(4,F)$ induced from the representation
$$
 \mat{A}{*}{}{cA'}\longmapsto\sigma(c)\pi(A)
$$
of $P(F)$. If $\chi$ is a character of $F^\times$ and $\pi$ is an admissible
representation of $\GSp(2,F)=\GL(2,F)$, then $\chi\rtimes\pi$ denotes the
representation of $\GSp(4,F)$ parabolically induced from the representation
$$
 \begin{bmatrix}x&*&*\\&A&*\\&&\det(A)x^{-1}\end{bmatrix}
 \longmapsto\chi(x)\pi(A)
$$
of $Q(F)$.

If $\Pi$ is an admissible representation of $\GSp(4,F)$ and
$\tau$ is a character of $F^\times$, then the \emph{twist} of $\Pi$ by $\tau$,
denoted $\tau\Pi$, is the representation $g\mapsto\tau(\lambda(g))\Pi(g)$,
where $\lambda$ is the multiplier homomorphism. The effect of twisting on
parabolically induced representations is as follows:
$$
 \tau(\chi_1\times\chi_2\rtimes\sigma)=\chi_1\times\chi_2\rtimes\tau\sigma,\quad
 \tau(\pi\rtimes\sigma)=\pi\rtimes\tau\sigma,\quad
 \tau(\chi\rtimes\pi)=\chi\rtimes\tau\pi.
$$
The non-supercuspidal, irreducible, admissible representations of $\GSp(4,F)$
have been classified by Sally and Tadi\'c in \cite{ST}. They determined
the irreducible subquotients of each representation parabolically
induced from an irreducible representation of $B$, $P$ or $Q$. In \cite{RS} this 
information was reorganized in the form of a table, which we reproduce
here as Table \ref{maintable}. The representations are organized in cases 
I - XI. Cases I - VI contain representations supported in $B$, cases VII - IX contain
those supported in $Q$, and cases X and XI contain representations
supported in $P$. For example, case I contains the irreducible, admissible 
representations of the form $\chi_1\times\chi_2\rtimes\sigma$. We refer to 
\cite[Sect.\ 2.2]{RS} for a precise description of the various cases.  

\subsection{Weil group representations}\label{weil}
We recall some basic facts about the Weil group $W_F$ and the Weil-Deligne group
$W'_F$ of $F$, referring to \cite{Roh} and \cite{T} for details.  
Recall from local Class Field Theory that the abelianized Weil group $W_F^{\rm ab}$ 
and $F^\times$ are isomorphic, which implies that the characters of $W_F$ and those 
of $F^\times$ can be identified. We will use the same symbol for a character of 
$F^\times$ and the corresponding character of $W_F$.  Representations of the 
Weil-Deligne group $W_F'$ are given by pairs $(\rho,N)$, where $\rho$ is a 
continuous homomorphism $W_F\rightarrow\GL(n,\C)$
and $N$ is a nilpotent complex $n\times n$ matrix for which
$$
 \rho(w)N\rho(w)^{-1}=\nu(w)N\qquad\text{for all }w\in W_F.
$$
If $\rho$ is a semisimple representation, then $(\rho,N)$ is called
\emph{admissible}. One attaches an $L$-factor $L(s,\varphi)$ to the 
pair  $\varphi=(\rho,N)$ as follows.  Let $\Phi\in W_F$ be an inverse Frobenius 
element and let $I={\rm Gal}(\bar F/F^{\rm un})\subset W_F$
be the inertia subgroup. Let $V_N=\ker(N)$,
$V^I=\{v\in V:\:\rho(g)v=v\;\text{for all }g\in I\}$ and $V^I_N=V^I\cap V_N$.
Then
\begin{equation}\label{LWDdefeq}
 L(s,\varphi)=\det\left(1-q^{-s}\rho(\Phi)\big|V^I_N\right)^{-1}.
\end{equation}
If $\varphi$ is a one-dimensional representation identified
with a character $\chi$ of $F^\times$, then
$$
 L(s,\varphi)=L(s,\chi)=\left\{\begin{array}{l@{\qquad\text{if }\chi}l}
 1&\text{ is ramified},\\ (1-\chi(\varpi)q^{-s})^{-1}
 &\text{ is unramified}.\end{array}\right.
$$
An \emph{$L$-parameter} for $\GSp(4,F)$ is essentially an equivalence class
of admissible homomorphisms $W'_F\rightarrow\GSp(4,\C)$; for the precise
definition see \cite[Sect.\ 2.4]{RS}. The conjectural local Langlands correspondence
assigns to each irreducible, admissible representation $\Pi$ of $\GSp(4,F)$ an
$L$-parameter $\varphi_\Pi$. It was shown in \cite[Sect.\ 2.4]{RS} that, for the
non-supercuspidal representations of $\GSp(4,F)$, there is a
unique way to make this assignment in such a way as to satisfy certain desirable properties
of the local Langlands correspondence. In what follows we shall always refer
to these unique parameters $\varphi_\Pi$ when we talk about the local Langlands
correspondence for the non-supercuspidal representations of $\GSp(4,F)$.
Their explicit forms are given in \cite[Sect.\ 2.4]{RS} and will be recalled below.

%%%%%%%%
\section{Computations of adjoint $L$-functions}
%%%%%%%%
We now go through the list of non-supercuspidal, irreducible, admissible
representations of $\GSp(4,F)$ and compute the adjoint $L$-functions 
of the $L$-parameters of these representations. 

\subsection{Cases supported in the minimal parabolic subgroup} \ 

\underline{Case I:} These are irreducible representations of the form
$\chi_1\times\chi_2\rtimes\sigma$, where $\chi_1$, $\chi_2$ and $\sigma$
are characters of $F^\times$. The condition for irreducibility is that
$\chi_1\neq\nu^{\pm1}$, $\chi_2\neq\nu^{\pm1}$ and $\chi_1\neq\nu^{\pm1}\chi_2^{\pm1}$.
The $L$-parameter of such a representation is given by the pair $(\rho,N)$,
where $N=0$ and
$$
 \rho(w)=\begin{bmatrix}(\chi_1\chi_2\sigma)(w)\\&(\chi_1\sigma)(w)\\
 &&(\chi_2\sigma)(w)\\&&&\sigma(w)\end{bmatrix}.
$$
The one-dimensional spaces spanned by the vectors in
(\ref{rootspaces1eq}) through (\ref{rootspaces4eq}) are preserved by the
action of $W_F$ on the $10$-dimensional space $\ssp(4)$ given by $\Ad_{\ssp(4)}\circ\rho$.
More precisely, $W_F$ acts on $L_\alpha$ by multiplication with $\alpha(\rho(w))$,
for each root $\alpha$. Furthermore, $W_F$ acts trivially on the diagonal torus
of $\ssp(4)$. Thus
\begin{align}\label{ILADeq}
 L(s,\chi_1\times\chi_2\rtimes\sigma,\Ad)=
 &L(s,1_{F^\times})^2L(s,\chi_1)L(s,\chi_1^{-1})L(s,\chi_2)L(s,\chi_2^{-1})\nonumber\\
 &L(s,\chi_1\chi_2)L(s,\chi_1^{-1}\chi_2^{-1})
  L(s,\chi_1\chi_2^{-1})L(s,\chi_1^{-1}\chi_2).
\end{align}

\underline{Case II:} Let $\chi$ and $\sigma$ be characters of $F^\times$ such that
$\chi^2\neq\nu^{\pm1}$ and $\chi\neq\nu^{\pm3/2}$. The induced representation
$\nu^{1/2}\chi\times\nu^{-1/2}\chi\rtimes\sigma$ has the two irreducible constituents
$\chi\St_{\GL(2)}\rtimes\sigma$ (type IIa) and $\chi\triv_{\GL(2)}\rtimes\sigma$
(type IIb). The $L$-parameter attached to $\chi\triv_{\GL(2)}\rtimes\sigma$
is $(\rho,N)$ with $N=0$ and
$$
 \rho(w)=\begin{bmatrix}(\chi^2\sigma)(w)\\&(\nu^{1/2}\chi\sigma)(w)\\
 &&(\nu^{-1/2}\chi\sigma)(w)\\&&&\sigma(w)\end{bmatrix}.
$$
Arguing similarly as in case I above, we obtain
\begin{align}\label{IIbLADeq}
 L(s,\chi\triv_{\GL(2)}\rtimes\sigma,\Ad)=
 &L(s,1_{F^\times})^2L(s,\chi^2)L(s,\chi^{-2})L(s,\nu)L(s,\nu^{-1})\nonumber\\
 &L(s,\chi\nu^{-1/2})L(s,\chi^{-1}\nu^{1/2})
  L(s,\chi\nu^{1/2})L(s,\chi^{-1}\nu^{-1/2}).
\end{align}

The $L$-parameter of the IIa type representation $\chi\St_{\GL(2)}\rtimes\sigma$
has the same semisimple part $\rho$, but $N=N_1$, where
\begin{equation}\label{N1defeq}
 N_1=\begin{bmatrix}0\\&0&1\\&&0\\&&&0\end{bmatrix}.
\end{equation}
Composing with the adjoint representation, the $10$-dimensional representation
of $W_F'$ whose $L$-factor we have to compute is $(\Ad_{\ssp(4)}\circ\rho,\,\ad(N_1))$.
To determine the $L$-factor we have to consider the restriction of
$\Ad_{\ssp(4)}\circ\rho$ to the kernel of $\ad(N_1)$; see (\ref{LWDdefeq}).
It is easy to see that
\begin{equation}\label{keradN1eq}
 \ker(\ad(N_1))=\langle\begin{bmatrix} 1\\&0\\&&0\\&&&-1\end{bmatrix},\,
 L_{2e_1},\,L_{e_1+e_2},\,L_{2e_2},\,L_{-e_1+e_2},\,L_{-2e_1}\rangle.
\end{equation}
The restriction of $\Ad_{\ssp(4)}\circ\rho$ to this $6$-dimensional space
decomposes in an obvious way into $1$-dimensional invariant subspaces, so
that the resulting $L$-factor is
\begin{align}\label{IIaLADeq}
 L(s,\chi\St_{\GL(2)}\rtimes\sigma,\Ad)=&L(s,1_{F^\times})L(s,\chi^2)
  L(s,\chi^{-2})\nonumber\\
 &L(s,\nu)L(s,\chi^{-1}\nu^{1/2})L(s,\chi\nu^{1/2}).
\end{align}

\underline{Case III:}
If $\chi$ and $\sigma$ are characters of $F^\times$ such that $\chi\neq1$ and
$\chi\neq\nu^{\pm2}$, then the induced representation $\chi\times\nu\rtimes
\nu^{-1/2}\sigma$ has two irreducible constituents $\chi\rtimes\sigma\St_{\GSp(2)}$
(type IIIa) and $\chi\rtimes\sigma\triv_{\GSp(2)}$ (type IIIb).
The $L$-parameter of $\chi\rtimes\sigma\triv_{\GSp(2)}$ is $(\rho,N)$ with
$N=0$ and
$$
 \rho(w)=\begin{bmatrix}(\nu^{1/2}\chi\sigma)(w)\\&(\nu^{-1/2}\chi\sigma)(w)\\
 &&(\nu^{1/2}\sigma)(w)\\&&&(\nu^{-1/2}\sigma)(w)\end{bmatrix}.
$$
Arguing as above, we find that
\begin{align}\label{IIIbLADeq}
 L(s,\chi\rtimes\sigma\triv_{\GSp(2)},\Ad)=
 &L(s,1_{F^\times})^2L(s,\chi)L(s,\chi^{-1})L(s,\nu)L(s,\nu^{-1})\nonumber\\
 &L(s,\chi\nu)L(s,\chi\nu^{-1})L(s,\chi^{-1}\nu)
  L(s,\chi^{-1}\nu^{-1}).
\end{align}
The $L$-parameter of $\chi\rtimes\sigma\St_{\GSp(2)}$ is $(\rho,N_4)$ with
the same $\rho$ and
\begin{equation}\label{N4defeq}
 N_4=\begin{bmatrix}0&1\\&0\\&&0&-1\\&&&0\end{bmatrix}.
\end{equation}
Composing with the adjoint representation, we obtain the representation of $W'_F$
given by $(\Ad_{\ssp(4)}\circ\rho,\ad(N_4))$.
It is easily computed that
\begin{equation}\label{keradN4eq}
 \ker(\ad(N_4))=\langle\begin{bmatrix} 1\\&1\\&&-1\\&&&-1\end{bmatrix},\,
 L_{-2e_2},\,L_{e_1-e_2},\,L_{2e_1}\rangle.
\end{equation}
Using the definition (\ref{LWDdefeq}) it follows that
\begin{align}\label{IIIaLADeq}
 L(s,\chi\St_{\GL(2)}\rtimes\sigma,\Ad)=&L(s,1_{F^\times})L(s,\nu)
  L(s,\nu\chi)L(s,\nu\chi^{-1}).
\end{align}

\underline{Case IV:}
Representations of type IV are the subquotients of $\nu^2\times\nu\rtimes\nu^{-3/2}\sigma$,
where $\sigma$ is a character of $F^\times$. The Langlands quotient is
$\sigma\triv_{\GSp(4)}$, a twist of the trivial representation (type IVd). Its $L$-parameter
is given by $(\rho,N)$ with $N=0$ and
$$
 \rho(w)=\begin{bmatrix}(\nu^{3/2}\sigma)(w)\\&(\nu^{1/2}\sigma)(w)\\
 &&(\nu^{-1/2}\sigma)(w)\\&&&(\nu^{-3/2}\sigma)(w)\end{bmatrix}.
$$
Arguing as before, we obtain
\begin{align}\label{IVdLADeq}
 L(s,\sigma\triv_{\GSp(4)},\Ad)=
 &L(s,1_{F^\times})^2L(s,\nu)^2L(s,\nu^{-1})^2L(s,\nu^2)L(s,\nu^{-2})\nonumber\\
 &L(s,\nu^3)L(s,\nu^{-3}).
\end{align}
The $L$-parameter
of the IVc type representation $L(\nu^{3/2}\St_{\GL(2)},\nu^{-3/2}\sigma)$ is
$(\rho,N_1)$ with $N_1$ as in (\ref{N1defeq}). It follows from (\ref{keradN1eq}) that
\begin{align}\label{IVcLADeq}
 L(s,L(\nu^{3/2}\St_{\GL(2)},\nu^{-3/2}\sigma),\Ad)=
 &L(s,1_{F^\times})L(s,\nu)L(s,\nu^{-1})L(s,\nu^2)\nonumber\\
 &L(s,\nu^3)L(s,\nu^{-3}).
\end{align}
The $L$-parameter
of the IVb type representation $L(\nu^2,\nu^{-1}\sigma\St_{\GSp(2)})$
is $(\rho,N_4)$ with $N_4$ as in (\ref{N4defeq}). 
It follows from (\ref{keradN4eq}) that
\begin{align}\label{IVbLADeq}
 L(s,L(\nu^2,\nu^{-1}\sigma\St_{\GSp(2)}),\Ad)=
 &L(s,1_{F^\times})L(s,\nu)L(s,\nu^{-1})L(s,\nu^3).
\end{align}
The $L$-parameter
of the IVa type representation $\sigma\St_{\GSp(4)}$ is $(\rho,N_5)$ with
\begin{equation}\label{N5defeq}
 N_5=\begin{bmatrix}0&1\\&0&1\\&&0&-1\\&&&0\end{bmatrix}.
\end{equation}
Easy computations show that
\begin{equation}\label{keradN5eq}
 \ker(\ad(N_5))=\langle L_{2e_1},\,L_{2e_2}+L_{e_1-e_2}\rangle.
\end{equation}
Thus
\begin{align}\label{IVaLADeq}
 L(s,\sigma\St_{\GSp(4)},\Ad)=&L(s,\nu)L(s,\nu^3).
\end{align}

\underline{Case V:} These are the irreducible subquotients of an induced representation
of the form $\nu\xi\times\xi\rtimes\nu^{-1/2}\sigma$, where $\xi$ is a non-trivial quadratic 
character of $F^\times$ and $\sigma$ is an arbitrary character of $F^\times$.
One of these subquotients is $L(\nu\xi,\xi\rtimes\nu^{-1/2}\sigma)$ (type IVd), and
its $L$-parameter is $(\rho,N)$ with $N=0$ and $\rho$ given by
$$
 \rho(w)=\begin{bmatrix}(\nu^{1/2}\sigma)(w)\\&(\nu^{1/2}\xi\sigma)(w)\\
 &&(\nu^{-1/2}\xi\sigma)(w)\\&&&(\nu^{-1/2}\sigma)(w)\end{bmatrix}.
$$
As in the other cases with $N=0$ one computes
\begin{align}\label{VdLADeq}
 L(s,L(\nu\xi,\xi\rtimes\nu^{-1/2}\sigma),\Ad)=
 &L(s,1_{F^\times})^2L(s,\nu)^2L(s,\nu^{-1})^2\nonumber\\
 &L(s,\xi)^2L(s,\nu\xi)L(s,\nu^{-1}\xi).
\end{align}
The $L$-parameter attached to the Vc type representation
$L(\nu^{1/2}\xi\St_{\GL(2)},\xi\nu^{-1/2}\sigma)$ is $(\rho,N_2)$ with the
same $\rho$ and
\begin{equation}\label{N2defeq}
 N_2=\begin{bmatrix}0&&&1\\&0\\&&0\\&&&0\end{bmatrix}.
\end{equation}
Computations show that
\begin{equation}\label{keradN2eq}
 \ker(\ad(N_2))=\langle\begin{bmatrix}0\\&1\\&&-1\\&&&0\end{bmatrix},\,
 L_{2e_1},\,L_{e_1+e_2},\,L_{2e_2},\,L_{-2e_2},\,L_{e_1-e_2}\rangle.
\end{equation}
Hence
\begin{align}\label{VcLADeq}
 L(s,L(\nu^{1/2}\xi\St_{\GL(2)},\xi\nu^{-1/2}\sigma),\Ad)=
 &L(s,1_{F^\times})L(s,\nu)^2L(s,\nu^{-1})\nonumber\\
 &L(s,\xi)L(s,\nu\xi).
\end{align}
The representation
$L(\nu^{1/2}\xi\St_{\GL(2)},\nu^{-1/2}\sigma)$ of type Vb is a $\xi$-twist of Vc. Since
adjoint $L$-functions are invariant under twists, its adjoint $L$-function
is the same as in (\ref{VcLADeq}).
The essentially square-integrable Va type representation
$\delta([\xi,\nu\xi],\nu^{-1/2}\sigma)$ has $L$-parameter $(\rho,N_3)$ with $\rho$
as before and
\begin{equation}\label{N3defeq}
 N_3=\begin{bmatrix}0&&&1\\&0&1\\&&0\\&&&0\end{bmatrix}.
\end{equation}
It is easy to compute that
\begin{equation}\label{keradN3eq}
 \ker(\ad(N_3))=\langle L_{2e_1},\,L_{e_1+e_2},\,L_{2e_2},\,L_{e_1-e_2}-L_{-e_1+e_2}\rangle.
\end{equation}
It follows that
\begin{align}\label{VaLADeq}
 L(s,\delta([\xi,\nu\xi],\nu^{-1/2}\sigma),\Ad)=
 &L(s,\nu)^2L(s,\xi)L(s,\nu\xi).
\end{align}

\underline{Case VI:} These are the irreducible subquotients of an induced representation
of the form $\nu\times1_{F^\times}\rtimes\nu^{-1/2}\sigma$, where
$\sigma$ is a character of $F^\times$. One of these irreducible subquotients
is the VId type representation $L(\nu,1_{F^\times}\rtimes\nu^{-1/2}\sigma)$.
Its $L$-parameter is $(\rho,N)$ with $N=0$ and
$$
 \rho(w)=\begin{bmatrix}(\nu^{1/2}\sigma)(w)\\&(\nu^{1/2}\sigma)(w)\\
 &&(\nu^{-1/2}\sigma)(w)\\&&&(\nu^{-1/2}\sigma)(w)\end{bmatrix}.
$$
The resulting adjoint $L$-function is
\begin{align}\label{VIdLADeq}
 L(s,L(\nu,1_{F^\times}\rtimes\nu^{-1/2}\sigma),\Ad)=
 &L(s,1_{F^\times})^4L(s,\nu)^3L(s,\nu^{-1})^3. 
\end{align}
The $L$-parameter of the
VIc type representation $L(\nu^{1/2}{\rm St}_{\GL(2)},\nu^{-1/2}\sigma)$
is $(\rho,N_1)$ with $N_1$ as in (\ref{N1defeq}). By (\ref{keradN1eq}),
\begin{align}\label{VIcLADeq}
 L(s,L(\nu^{1/2}{\rm St}_{\GL(2)},\nu^{-1/2}\sigma),\Ad)=
 &L(s,1_{F^\times})^2 L(s,\nu)^3 L(s,\nu^{-1}). 
\end{align}
The remaining irreducible subquotients are
the generic $\tau(S,\nu^{-1/2}\sigma)$ and the non-generic $\tau(T,\nu^{-1/2}\sigma)$.
Both of these are tempered representations and they constitute an $L$-packet.
Their common $L$-parameter is $(\rho,N_3)$ with $\rho$ as above and
$N_3$ as in (\ref{N3defeq}). By (\ref{keradN3eq}),
\begin{align}\label{VIabLADeq}
 L(s,\tau(S/T,\nu^{-1/2}\sigma),\Ad)=
 &L(s,1_{F^\times}) L(s,\nu)^3.
\end{align}

\subsection{Cases supported in the Klingen parabolic subgroup} \ 

\underline{Case VII:} These representations are the irreducible
representations of the form $\chi\rtimes\pi$, where $\chi$ is a character of $F^\times$
and $\pi$ is a supercuspidal irreducible admissible representation of $\GL(2,F)$.
If $\mu:\:W_F\rightarrow\GL(2,\C)$ is the $L$-parameter of $\pi$, then
$\chi\rtimes\pi$ has $L$-parameter $(\rho,N)$ with $N=0$ and
\begin{equation}\label{VIILparametereq}
 \rho(w)=\mat{\chi(w)\det(\mu(w))\mu(w)'}{}{}{\mu(w)}\in\GSp(4,\C).
\end{equation}
To compute the adjoint $L$-function of this parameter,
we identify the Siegel Levi $M_P$ in $\GSp(4,\C)=\widehat{\GSp}(4,F)$ with $\GL(2,\C)\times\GL(1,\C)$ via
\begin{equation}\label{MPidenteq}
 (A,x)\longmapsto\mat{xA'}{}{}{A}\qquad(A\in\GL(2,\C),\;x\in\GL(1,\C)).
\end{equation}
We have to decompose the Lie algebra $\mathfrak{sp}(4)$ into irreducible
representations of $M_P$. It is easy to see that
\begin{align}\label{sp4Pdecompeq}
 \nonumber\mathfrak{sp}(4)&=
 \underbrace{\C\begin{bmatrix}1\\&1\\&&-1\\&&&-1\end{bmatrix}}_{\text{invariant}}
 \oplus\underbrace{\C L_{-e_1+e_2}\oplus\C\begin{bmatrix}1\\&-1\\&&1\\&&&-1\end{bmatrix}
  \oplus\C L_{e_1-e_2}}_{\text{invariant}}\\[1ex]
 &\oplus\underbrace{\C L_{2e_2}\oplus\C L_{e_1+e_2}\oplus\C L_{2e_1}}_{\text{invariant}}
  \quad\oplus\quad
 \underbrace{\C L_{-2e_1}\oplus\C L_{-e_1-e_2}\oplus\C L_{-2e_2}}_{\text{invariant}}.
\end{align}
The representation on the $1$-dimensional invariant subspace is the trivial
representation.
The representation on $\C L_{2e_2}\oplus\C L_{e_1+e_2}\oplus\C L_{2e_1}$ is
$$
 \underbrace{({\rm det}^{-2}\otimes{\rm Sym}^2)}_{\text{representation of }\GL(2,\C)}
  \otimes\;\std_{\GL(1)}.
$$
The representation on $\C L_{-e_1+e_2}\oplus\C\begin{bmatrix}1\\&-1\\&&1\\&&&-1\end{bmatrix}
\oplus\C L_{e_1-e_2}$ is
$$
 ({\rm det}^{-1}\otimes{\rm Sym}^2)\otimes\;{\rm triv}_{\GL(1)}.
$$
The representation on $\C L_{-2e_1}\oplus\C L_{-e_1-e_2}\oplus\C L_{-2e_2}$ is
$$
 {\rm Sym}^2\otimes\;{\rm std}_{\GL(1)}^{-1}.
$$
Using ${\rm Sym}^2=\det\otimes{\rm Ad}_{\GL(2)}$ as representations of $\GL(2,\C)$, 
we can rewrite these three-dimensional
representations as
\begin{align*}
 ({\rm det}^{-1}\otimes{\rm Ad}_{\GL(2)})&\otimes\;{\rm std}_{\GL(1)},\\
 {\rm Ad_{\GL(2)}}&\otimes\;{\rm triv}_{\GL(1)},\\
 ({\rm det}\otimes{\rm Ad}_{\GL(2)})&\otimes\;{\rm std}_{\GL(1)}^{-1}.
\end{align*}
Via the identification (\ref{MPidenteq}), we consider $\rho$ as a homomorphism
$W_F\rightarrow\GL(2,\C)\times\GL(1,\C)$. As such we have
$\rho=\mu\times\chi\omega_\pi$; note that $\det\circ\mu=\omega_\pi$.
For the resulting $L$-functions we have the following lemma.
\begin{lemma}\label{L2Lrellemma}
 For a character $\chi$ of $F^\times$ and an irreducible admissible representation
 $\pi$ of $\GL(2,F)$, let
 $$
  L_2(s,\pi,\chi)=\frac{L(s,(\chi\pi)\times\tilde\pi)}{L(s,\chi)},
 $$
 as in \cite{GJ}. Then
 $$
  L(s,(({\rm det}^{-1}\otimes{\rm Ad}_{\GL(2)})\otimes\std_{\GL(1)})
 \circ(\mu\times(\chi\omega_\pi)))=L_2(s,\pi,\chi).
 $$
 and
 $$
  L(s,(({\rm det}\otimes{\rm Ad}_{\GL(2)})\otimes\std_{\GL(1)}^{-1})
 \circ(\mu\times(\chi\omega_\pi)))=L_2(s,\pi,\chi^{-1}).
 $$
 Here, $\mu:\:W'_F\rightarrow\GL(2,\C)$ is the $L$-parameter of $\pi$.
\end{lemma}
\begin{proof}
We have
$$
 \std_{\GL(2)}\otimes\std_{\GL(2)}=\det\otimes({\rm Ad}_{\GL(2)}\oplus\triv_{\GL(2)}),
$$
and hence
\begin{align*}
 L(s,(\chi\pi)\times\tilde\pi))&=L(s,(\chi\omega_\pi^{-1})\pi\times\pi)\\
 &=L(s,\chi\cdot({\rm det}^{-1}\circ\mu)(\mu\otimes\mu))\\
 &=L(s,\chi\cdot({\rm Ad}_{\GL(2)}\oplus\triv_{\GL(2)})\circ\mu)\\
 &=L(s,\chi)L(s,\chi\cdot({\rm Ad}_{\GL(2)}\circ\mu))\\
 &=L(s,\chi)L(s,(\chi\omega_\pi)\cdot(({\rm det}^{-1}{\rm Ad}_{\GL(2)})\circ\mu))\\
 &=L(s,\chi)L(s,(({\rm det}^{-1}{\rm Ad}_{\GL(2)})\otimes\std_{\GL(1)})
 \circ(\mu\times(\chi\omega_\pi))).
\end{align*}
\end{proof}
\begin{rem}
One can write the $L$-function $L_2$ in the standard notation of 
Langlands $L$-functions as 
$$ 
L_2(s,\pi,\chi) = L(s,\pi,{\rm Sym}^2\otimes (\omega_\pi^{-1}\chi)) = L(s,\pi,{\rm Ad}_{\GL(2)}\otimes\chi), 
$$ 
where $\omega_\pi$ denotes the central character of $\pi$ and we use the same symbol for both 
characters of $F^\times$ and the corresponding characters of $W_F$ as in \ref{weil}. 
This means that $L_2$ is a twisted symmetric squre, or equivalently, a twisted adjoint $L$-function 
of $\GL(2)$. (The adjoint $L$-function is sometimes also referred to as the adjoint square $L$-function.) 
We will use the latter notation in our final formulas below.
\end{rem}

It follows that
\begin{align}\label{VIILADeq}
\nonumber 
 L(s,\chi\rtimes\pi,\Ad)=
 &L(s,1_{F^\times})L(s,\pi,{\rm Ad}_{\GL(2)}) \cdot \\
 &L(s,\pi,{\rm Ad}_{\GL(2)}\otimes\chi)L(s,\pi,{\rm Ad}_{\GL(2)}\otimes\chi^{-1}).
\end{align} 

\underline{Case VIII:} If $\pi$ is a supercuspidal irreducible admissible
representation of $\GL(2,F)$, then the induced representation $1_{F^\times}\rtimes\pi$
is a direct sum of two irreducible constituents $\tau(S,\pi)$ (type VIIIa) and
$\tau(T,\pi)$ (type VIIIb). Both irreducible constituents are tempered, but only
VIIIa is generic. These two representations constitute an $L$-packet. Their
common $L$-parameter is $(\rho,N)$ with $N=0$ and
$$
 \rho(w)=\mat{\det(\mu(w))\mu(w)'}
 {}{}{\mu(w)}\in\GSp(4,\C). 
$$
Here, $\mu:\:W_F\rightarrow\GL(2,\C)$ is the parameter of $\pi$. The calculation
of the adjoint $L$-function of this parameter is exactly as in Case VII. 
The result is
\begin{align}\label{VIIILADeq}
 L(s,1_{F^\times}\rtimes\pi,\Ad)=
 &L(s,1_{F^\times})L(s,\pi,{\rm Ad}_{\GL(2)})^3.
\end{align} 

\underline{Case IX:} These are the irreducible constituents of induced representations
of the form $\nu\xi\rtimes\nu^{-1/2}\pi$, where $\xi$ is a non-trivial quadratic
character of $F^\times$, and where $\pi$ is a supercuspidal representation of
$\GL(2,F)$ for which $\xi\pi=\pi$. The generic constituent is denoted by
$\delta(\nu\xi,\nu^{-1/2}\pi)$ (type IXa), and the non-generic constituent is
denoted by $L(\nu\xi,\nu^{-1/2}\pi)$ (type IXb). The $L$-parameter of
$L(\nu\xi,\nu^{-1/2}\pi)$ is $(\rho,N)$, where $N=0$ and
\begin{equation}\label{LparameterIXeq}
 \rho(w)=\mat{\xi(w)\nu^{1/2}(w)\det(\mu(w))\mu'(w)}
 {}{}{\nu^{-1/2}(w)\mu(w)}. 
\end{equation}
Here, $\mu:\:W_F\rightarrow\GL(2,\C)$ is the $L$-parameter of $\pi$.
The computation of the adjoint $L$-function of this representation is very similar
to type VII above. The result is
\begin{align}\label{IXbLADeq}
 L(s,L(\nu\xi,\nu^{-1/2}\pi),\Ad)=
 &L(s,1_{F^\times})L(s,\nu^{-1/2}\pi,{\rm Ad}_{\GL(2)}) \nonumber\\
 &\cdot L_2(s,\nu^{-1/2}\pi,\xi\nu)L_2(s,\nu^{-1/2}\pi,\xi\nu^{-1})\nonumber\\
= &L(s,1_{F^\times})L(s,\pi,{\rm Ad}_{\GL(2)}) \nonumber\\
 &\cdot L(s,\pi,{\rm Ad}_{\GL(2)}\otimes\xi\nu)L(s,\pi,{\rm Ad}_{\GL(2)}\otimes\xi\nu^{-1}).
\end{align}

The $L$-parameter of $\delta(\nu\xi,\nu^{-1/2}\pi)$ is $(\rho,N)$, where $\rho$
is as above and $N$ is defined as follows. By \cite[Lemma 2.4.1]{RS} there
exists a symmetric invertible matrix $S\in\GL(2,\C)$ such that
\begin{equation}\label{RS241eq}
 ^t\mu(w)S\mu(w)=\xi(w)\det(\mu(w))S\qquad\text{for all }w\in W_F.
\end{equation}
Then $N=\mat{0}{B}{0}{0}$ with $B=\mat{0}{1}{1}{0}S$. We have to consider the
action of $W_F$ on $\ker({\rm ad}(N))$ via ${\rm Ad}\circ\rho$.
It is clear that $\ker({\rm ad}(N))$ contains the subspace
$\C L_{2e_2}\oplus\C L_{e_1+e_2}\oplus\C L_{2e_1}$ appearing in
(\ref{sp4Pdecompeq}). The operator ${\rm ad}(N)$ induces a linear map
\begin{equation}\label{splinearmap1eq}
 \ssp(4)\supset\mat{*}{0}{0}{*}\longrightarrow\mat{0}{*}{0}{0}\subset\ssp(4).
\end{equation}
The domain of this linear map is $4$-dimensional, and the target space is
$3$-dimensional. It is easy to see that, since $S$ is invertible, this
linear map is surjective. It follows that there exists a non-zero matrix
$A_0\in M(2\times2,\C)$, unique up to scalars, for which
\begin{equation}\label{A0defeq}
 A_0B=-B\mat{0}{1}{1}{0}\,^t\!A_0\mat{0}{1}{1}{0}.
\end{equation}
In fact, a calculation verifies that
\begin{equation}\label{A0Seq}
 A_0=\mat{}{1}{1}{}S\mat{-1}{}{}{1}
\end{equation}
is such a matrix.
Furthermore, we get $\dim(\ker({\rm ad}(N)))\geq4$ and
$\dim({\rm im}({\rm ad}(N)))\geq3$. In fact, we claim that
$$
 \dim(\ker({\rm ad}(N)))=4\qquad\text{and}\qquad
 \dim({\rm im}({\rm ad}(N)))=6.
$$
By what we already proved, it is enough to show that
$\dim({\rm im}({\rm ad}(N)))\geq6$. It is easy to see that ${\rm ad}(N)$
induces an injective linear map
\begin{equation}\label{splinearmap2eq}
 \ssp(4)\supset\mat{0}{0}{\ast}{0}\longrightarrow\mat{*}{0}{0}{*}\subset\ssp(4).
\end{equation}
It follows that the intersection of ${\rm im}({\rm ad}(N))$ with the Siegel Levi
is at least $3$-dimensional. Since ${\rm im}({\rm ad}(N))$
also contains the image of the map (\ref{splinearmap1eq}), it follows that
we have indeed $\dim({\rm im}({\rm ad}(N)))\geq6$. This proves our claim.
We showed that
$$
 \ker({\rm ad}(N))=\langle L_{2e_1},L_{e_1+e_2},L_{2e_2},\mat{A_0}{}{}{A_0'}\rangle,
 \qquad A_0'=-\mat{0}{1}{1}{0}\,^t\!A_0\mat{0}{1}{1}{0}.
$$
The action of $W_F$ preserves the Siegel Levi of $\ssp(4)$, and therefore
the one-dimensional space spanned by $\mat{A_0}{}{}{A_0'}$. Hence, 
$$
 \rho(w)\mat{A_0}{}{}{A_0'}\rho(w)^{-1}=\eta(w)\mat{A_0}{}{}{A_0'}
$$
for a character $\eta$ of $W_F$. In fact, using (\ref{RS241eq}), it is easy to
see that $\eta=\xi$.
This one-dimensional subspace therefore contributes a factor $L(s,\xi)$ to the $L$-function.
The $L$-factor resulting from the action of $W_F$ on
$\C L_{2e_2}\oplus\C L_{e_1+e_2}\oplus\C L_{2e_1}$ has been computed before;
see Lemma \ref{L2Lrellemma}. We finally get
\begin{align}\label{IXaLADeq}
 L(s,\delta(\nu\xi,\nu^{-1/2}\pi),\Ad)
 &=L(s,\xi)L_2(s,\nu^{-1/2}\pi,\xi\nu)\nonumber\\
 &=L(s,\xi)L(s,\pi,{\rm Ad}_{\GL(2)}\otimes\xi\nu). 
\end{align}

\subsection{Cases supported in the Siegel parabolic subgroup} \ 

\underline{Case X:} This case consists of the irreducible admissible
representations of $\GSp(4,F)$ of the form $\pi\rtimes\sigma$, where $\pi$
is a supercuspidal, irreducible representation of $\GL(2,F)$ and $\sigma$ is a character
of $F^\times$. The condition for irreducibility is that the central character
$\omega_\pi$ of $\pi$ is not equal to $\nu^{\pm1}$. If $\mu:\:W_F\rightarrow\GL(2,\C)$
is the $L$-parameter of $\pi$, then the $L$-parameter of $\pi\rtimes\sigma$ is
$(\rho,N)$ with $N=0$ and
\begin{align}\label{localparameterXeq}
 \rho(w)=\begin{bmatrix}\sigma(w)\det(\mu(w))&&\\
 &\sigma(w)\mu(w)\\&&\sigma(w)\end{bmatrix}.
\end{align}
In particular, the image of $\rho$ is contained in $M_Q$,
the standard Levi subgroup of the Klingen parabolic.
It is easy to see that the restriction of the adjoint representation of $\GSp(4,\C)$
to $M_Q$ decomposes into the following invariant subspaces:
\begin{align}\label{sp4Qdecompeq}
 \ssp(4)&=\underbrace{\C\begin{bmatrix}1\\&0\\&&0\\&&&-1\end{bmatrix}}_{\text{invariant}}
 \quad\oplus\quad
 \underbrace{\C L_{2e_2}\oplus\C\begin{bmatrix}0\\&1\\&&-1\\&&&0\end{bmatrix}
 \oplus\C L_{-2e_2}}_{\text{invariant}}\nonumber\\
 &\oplus\quad\underbrace{\C L_{e_1+e_2}\oplus\C L_{e_1-e_2}}_{\text{invariant}}
  \quad\oplus\quad
  \underbrace{\C L_{-e_1+e_2}\oplus\C L_{-e_1-e_2}}_{\text{invariant}}\nonumber\\
 &\oplus\quad\underbrace{\C L_{2e_1}}_{\text{invariant}}\quad\oplus\quad
  \underbrace{\C L_{-2e_1}}_{\text{invariant}}.
\end{align}
The action of $W_F$ via ${\rm Ad}\circ\rho$ on the first invariant subspace is
trivial. The action on the second invariant subspace is ${\rm Ad}_{\sl(2)}\circ\mu$.
The action on the third invariant subspace is $\std_{\GL(2)}\circ\mu$.
The action on the fourth invariant subspace is the twist of the previous one
by $\det\circ\mu^{-1}$.
And the action on the last two invariant subspaces is via 
$\det\circ\mu$ and its inverse, respectively. Hence we get
$$
 L(s,\pi\rtimes\sigma,\Ad)=
 L(s,1_{F^\times})L(s,\pi,{\rm Ad}_{\GL(2)})L(s,\pi)L(s,\omega_\pi^{-1}\pi)L(s,\omega_\pi)
 L(s,\omega_\pi^{-1}).
$$
Since $\pi$ is supercuspidal, $L(s,\pi)=L(s,\omega_\pi^{-1}\pi)=1$, so that
\begin{align}\label{XLADeq}
 L(s,\pi\rtimes\sigma,\Ad)=
 &L(s,1_{F^\times})L(s,\pi,{\rm Ad}_{\GL(2)})L(s,\omega_\pi)L(s,\omega_\pi^{-1}).
\end{align}

\underline{Case XI:}
Let $\pi$ be a supercuspidal representation of $\GL(2,F)$ with $\omega_\pi=1$
and $\sigma$ a character of $F^\times$. Then $\nu^{1/2}\pi\rtimes\nu^{-1/2}\sigma$
decomposes into the XIa type representation $\delta(\nu^{1/2}\pi,\nu^{-1/2}\sigma)$
and the XIb type representation $L(\nu^{1/2}\pi,\nu^{-1/2}\sigma)$.
The Langlands quotient $L(\nu^{1/2}\pi,\nu^{-1/2}\sigma)$ has $L$-parameter
$(\rho,N)$ with $N=0$ and
\begin{align}\label{localparameterXIeq}
 \rho(w)=\begin{bmatrix}\sigma(w)\nu^{1/2}(w)&&\\
 &\sigma(w)\mu(w)\\&&\sigma(w)\nu^{-1/2}(w)\end{bmatrix}.
\end{align}
The computation of the adjoint $L$-function is the same as in Case X. The result is
\begin{align}\label{XIbLADeq}
 L(s,L(\nu^{1/2}\pi,\nu^{-1/2}\sigma),\Ad)=
 &L(s,1_{F^\times})L(s,\pi,{\rm Ad}_{\GL(2)})L(s,\nu)L(s,\nu^{-1}).
\end{align} 

The $L$-parameter of the XIa type representation
$\delta(\nu^{1/2}\pi,\nu^{-1/2}\sigma)$ is $(\rho,N_2)$ with the same $\rho$
as above and $N_2$ as defined in (\ref{N2defeq}). By (\ref{keradN2eq}),
we have to consider the restriction of ${\rm Ad}\circ\rho$ to the second,
third and fourth invariant subspace in (\ref{sp4Qdecompeq}). It follows that
\begin{align}\label{XIaLADeq}
 L(s,\delta(\nu^{1/2}\pi,\nu^{-1/2}\sigma),\Ad)=
 &L(s,\nu)L(s,\pi,{\rm Ad}_{\GL(2)}).
\end{align}
%%%%%%%%
\section{Generic Criterion}\label{generic}
%%%%%%%%
As a corollary of our computations we prove Theorem \ref{GP-R} below, which is a special case of a 
conjecture of Gross-Prasad and Rallis for non-supercuspidal representations of $\GSp(4,F)$. 

\begin{lemma}\label{L2poleslemma}
 Let $\pi$ be a supercuspidal representation
 of $\GL(2,F)$. Then the $L$-function $L_2(s,\pi,\chi)=L(s,\pi,{\rm Ad}_{\GL(2)}\otimes\chi)$ 
 in Lemma \ref{L2Lrellemma} has a pole at $s=1$
 if and only if $\chi=\nu^{-1}\xi$ with $\xi$ a non-trivial quadratic character for which 
 $\xi\pi\cong\pi$. In case of a pole, that pole must be simple.
\end{lemma}

\begin{proof} Assume that 
 $$
  L_2(s,\pi,\chi)=\frac{L(s,(\chi\pi)\times\tilde\pi)}{L(s,\chi)}
 $$
has a pole at $s=1$. Then
$L(s,(\chi\pi)\times\tilde\pi)$ has a pole at $s=1$. By \cite[Prop.\ (1.2)]{GJ} this 
implies that $\nu\chi\pi\cong\pi$. Taking central characters shows that 
$\chi=\nu^{-1}\xi$ with a quadratic character $\xi$. Since the pole of 
$L(s,(\chi\pi)\times\tilde\pi)$ is simple by \cite[Prop.\ (1.2)]{GJ}, our hypothesis 
implies that the function $L(s,\chi)$ cannot have a pole at $s=1$. Hence 
$\xi$ is non-trivial.  

Conversely, if $\chi=\nu^{-1}\xi$ with $\xi$ a non-trivial quadratic character 
for which $\xi\pi\cong\pi$, then $L_2(s,\pi,\chi)$ has a simple pole 
at $s=1$ by \cite[Prop.\ (1.2)]{GJ}.
\end{proof} 

\begin{theorem}\label{GP-R} 
Let $\varphi$ be the $L$-parameter of a non-supercuspidal, irreducible, admissible
representation of $\GSp(4,F)$ as above. Then
the $L$-function $L(s,\varphi,{\rm Ad})$ is holomorphic at $s=1$ if and only 
if one of the $L$-indistinguishable representations with $L$-parameter 
$\varphi$ listed in Table \ref{maintable} is generic.  
\end{theorem} 

\begin{proof} 
Among the representations listed in Table \ref{maintable} in each group the top one
(type ``a'') is generic. 
We verify that their adjoint $L$-functions are holomorphic at $s=1$ while the adjoint 
$L$-function of all the other representations do indeed have poles at $s=1$.  
Recall that the local factor $L(s,\chi)$ is always non-zero and it has a pole at 
$s=1$ if and only if $\chi=\nu^{-1}$.  We now go through the list and determine 
the order of the possible pole at $s=1$ using the irreducibility conditions for 
each case. The results are summarized in Table \ref{polestable}. 

In case I the irreducibility conditions $\chi_1\neq\nu^{\pm1}$, 
$\chi_2\neq\nu^{\pm1}$ and $\chi_1\neq\nu^{\pm1}\chi_2^{\pm1}$ 
imply that the $L$-function (\ref{ILADeq}) has no pole at $s=1$.

In case IIb, the factor $L(s,\nu^{-1})$ in (\ref{IIbLADeq}) contributes a simple pole at $s=1$ and 
the conditions $\chi^2\neq\nu^{\pm1}$ and $\chi\neq\nu^{\pm3/2}$ imply that none of the other 
factors contributes a pole at $s=1$. Also, it follows from (\ref{IIaLADeq}) that the adjoint 
$L$-function of a generic representation of type IIa has no pole at $s=1$.

The $L$-function in (\ref{IIIbLADeq}) for case IIIb has a double pole at $s=1$ if $\chi=\nu^{\pm1}$, and a simple pole otherwise.
Since $\chi\neq\nu^{\pm2}$, it follows from (\ref{IIIaLADeq}) that the 
adjoint $L$-function of a generic representation of type IIIa has no pole at $s=1$.

The adjoint $L$-function for cases IVa--IVd are, respectively, given in 
(\ref{IVaLADeq}), (\ref{IVbLADeq}), (\ref{IVcLADeq}), and (\ref{IVdLADeq}).
Clearly, the first has no pole, the second and third have simple poles,
and the fourth has a double pole at $s=1$. 

Similarly the adjoint $L$-function for case Va is given in (\ref{VaLADeq}), for cases 
Vb and Vc in  (\ref{VcLADeq}), and for case Vd in (\ref{VdLADeq}). Again, 
the first has no pole, the second and third have a simple pole, and  
the last a double pole at $s=1$. 

The representations in VIa and VIb are in the same $L$-packet. Their adjoint 
$L$-function, given in (\ref{VIabLADeq}), is holomorphic at $s=1$. On the 
other hand, the adjoint $L$-function of VIc is given in (\ref{VIcLADeq}) and has 
a simple pole at $s=1$. The adjoint $L$-function of VId is given
in (\ref{VIdLADeq}) with a triple pole at $s=1$. 

For case VII note that if we had $\chi=\nu^{-1}\xi$ with a non-trivial quadratic 
character $\xi$ for which $\xi\pi\cong\pi$, then $\chi\rtimes\pi$ would reduce 
and would therefore not be of type VII, but of type IX. Therefore, Lemma \ref{L2poleslemma}
implies that $L(s,\chi\rtimes\pi,\Ad)$, given by (\ref{VIILADeq}), has no pole at $s=1$ 
(note that $L(s,\pi,{\rm Ad}_{\GL(2)})$ is holomorphic at $s=1$ since $\pi$ is generic).

Cases VIIIa and VIIIb constitute an $L$-packet with VIIIa generic. Their adjoint 
$L$-function, given by (\ref{VIIILADeq}), is holomorphic at $s=1$ by 
Lemma \ref{L2poleslemma}.

Case IXb has the adjoint $L$-function given in (\ref{IXbLADeq}). By Lemma 
\ref{L2poleslemma} this $L$-function has a simple pole at $s=1$, coming from 
the factor $L(s,\pi,{\rm Ad}_{\GL(2)}\otimes\xi\nu^{-1})$. The adjoint $L$-function of case IXa 
is given in (\ref{IXaLADeq}). By Lemma \ref{L2poleslemma} the factor 
$L(s,\pi,{\rm Ad}_{\GL(2)}\otimes\xi\nu)$, and therefore $L(s,\delta(\nu\xi,\nu^{-1/2}\pi),\Ad)$,
has no pole at $s=1$.

The adjoint $L$-function for case X is given in (\ref{XLADeq}). Since
$\omega_\pi\neq\nu^{\pm1}$, this function is holomorphic at $s=1$. 

Finally, the adjoint $L$-functions for cases XIa and XIb are given in 
(\ref{XIaLADeq}) and (\ref{XIbLADeq}), respectively. The former is holomorphic 
at $s=1$ while the latter has a simple pole there. 
\end{proof}

\begin{rem} Cases Va and XIa are expected to have non-generic supercuspidal 
representations in their $L$-packets. Also, cases VIa and VIb as well as VIIIa and VIIIb 
constitute $L$-packets.  $L$-packets of all the other representations in Table \ref{maintable} 
are singletons.  
\end{rem}

%%%%%%%%
\newpage

%\begin{appendix}
\begin{table}
\centering
\caption{Non-supercuspidal representations of $\GSp(4,F)$}
\label{maintable}
\vspace{-3ex}
$$
\renewcommand{\arraystretch}{1.5}
 \begin{array}{|c|c|c|c|c|c|}
  \hline&&\mbox{constituent of}&\mbox{representation}
   &\mbox{centr.\ char.} 
   &\mbox{generic}\\\hline\hline
% I
  {\rm I}&&\multicolumn{2}{|c|}{\chi_1\times\chi_2\rtimes\sigma\quad
   \mbox{(irreducible)}}
   &\chi_1\chi_2\sigma^2  
   &\bullet\\\hline
% II
  &\mbox{a}&\nu^{1/2}\chi\times\nu^{-1/2}\chi\rtimes\sigma
   &\chi\St_{\GL(2)}\rtimes\sigma
   &  
   &\bullet\\\cline{4-4}\cline{6-6}
  \raisebox{2.5ex}[-2.5ex]{II}
   &\mbox{b}&(\chi^2\neq\nu^{\pm1},\chi\neq\nu^{\pm3/2})
   &\chi\triv_{\GL(2)}\rtimes\sigma
   & \raisebox{2.5ex}[-2.5ex]{$\chi^2\sigma^2$}   
   &\\\hline
% III
  &\mbox{a}&\chi\times\nu\rtimes\nu^{-1/2}\sigma
   &\chi\rtimes\sigma\St_{\GSp(2)}
   &  
   &\bullet \\\cline{4-4}\cline{6-6}
  \raisebox{2.5ex}[-2.5ex]{III}
  &\mbox{b}&(\chi\notin\{1,\nu^{\pm2}\})
   &\chi\rtimes\sigma\triv_{\GSp(2)}
   & \raisebox{2.5ex}[-2.5ex]{$\chi\sigma^2$} 
   &\\\hline
% IV
  &\mbox{a}&&\sigma\St_{\GSp(4)}
  &  
  &\bullet\\\cline{4-4}\cline{6-6}
  &\mbox{b}&&L(\nu^2,\nu^{-1}\sigma\St_{\GSp(2)})
   &  
   &\\ \cline{4-4}\cline{6-6}
  \raisebox{2.5ex}[-2.5ex]{IV}&
   \mbox{c}&\raisebox{2.5ex}[-2.5ex]{$\nu^2\times\nu\rtimes\nu^{-3/2}\sigma$}
   &L(\nu^{3/2}\St_{\GL(2)},\nu^{-3/2}\sigma)
   &   \raisebox{2.5ex}[-2.5ex]{$\sigma^2$} 
   &\\\cline{4-4}\cline{6-6}
  &\mbox{d}&&\sigma\triv_{\GSp(4)}
   &  
   &\\\hline
% V
  &\mbox{a}&&\delta([\xi ,\nu\xi ],\nu^{-1/2}\sigma)
  &  
  &\bullet\\\cline{4-4}\cline{6-6}
  &\mbox{b}&\nu\xi \times\xi \rtimes\nu^{-1/2}\sigma&
   L(\nu^{1/2}\xi \St_{\GL(2)},\nu^{-1/2}\sigma)
   &  
   &\\\cline{4-4}\cline{6-6}
  \raisebox{2.5ex}[-2.5ex]{V}&\mbox{c}&(\xi ^2=1,\:\xi \neq1)
   &L(\nu^{1/2}\xi \St_{\GL(2)},\xi \nu^{-1/2}\sigma)
   & \raisebox{2.5ex}[-2.5ex]{$\sigma^2$} 
   &\\\cline{4-4}\cline{6-6}
  &\mbox{d}&&L(\nu\xi ,\xi \rtimes\nu^{-1/2}\sigma)
   &  
   &\\\hline
% VI
  &\mbox{a}&&\tau(S,\nu^{-1/2}\sigma)
   &  
   &\bullet\\\cline{4-4}\cline{6-6}
  &\mbox{b}&&\tau(T,\nu^{-1/2}\sigma)
   &  
   &\\\cline{4-4}\cline{6-6}
  \raisebox{2.5ex}[-2.5ex]{VI}&\mbox{c}
   &\raisebox{2.5ex}[-2.5ex]{$\nu\times1_{F^\times}\rtimes\nu^{-1/2}\sigma$}
   &L(\nu^{1/2}\St_{\GL(2)},\nu^{-1/2}\sigma)
   & \raisebox{2.5ex}[-2.5ex]{$\sigma^2$} 
   &\\\cline{4-4}\cline{6-6}
  &\mbox{d}&&L(\nu,1_{F^\times}\rtimes\nu^{-1/2}\sigma)
   &  
   &\\\hline\hline
% VII 
  {\rm VII}&&\multicolumn{2}{|c|}{\chi\rtimes\pi\quad
   \mbox{(irreducible)}}
   & \chi\omega_\pi 
   &\bullet\\\hline
% VIII
  &\mbox{a}&
   &\tau(S,\pi)
   &  
   &\bullet\\\cline{4-4}\cline{6-6}
  \raisebox{2.5ex}[-2.5ex]{VIII}
  &\mbox{b}&\raisebox{2.5ex}[-2.5ex]{$1_{F^\times}\rtimes\pi$}
   &\tau(T,\pi)
   & \raisebox{2.5ex}[-2.5ex]{$\omega_\pi$} 
   &\\\hline
% IX
  &\mbox{a}&\nu\xi \rtimes\nu^{-1/2}\pi
   &\delta(\nu\xi ,\nu^{-1/2}\pi)
   &  
   &\bullet\\\cline{4-4}\cline{6-6}
  \raisebox{2.5ex}[-2.5ex]{IX}
  &\mbox{b}&(\xi \neq1,\:\xi \pi=\pi)
   &L(\nu\xi ,\nu^{-1/2}\pi)
   & \raisebox{2.5ex}[-2.5ex]{$\omega_\pi\xi$} 
   &\\\hline\hline
% X 
  {\rm X}&&\multicolumn{2}{|c|}{\pi\rtimes\sigma\quad
   \mbox{(irreducible)}}
   & \omega_\pi\sigma^2 
   &\bullet\\\hline
% XI
  &\mbox{a}&\nu^{1/2}\pi\rtimes\nu^{-1/2}\sigma
   &\delta(\nu^{1/2}\pi,\nu^{-1/2}\sigma)
   &  
   &\bullet\\\cline{4-4}\cline{6-6}
  \raisebox{2.5ex}[-2.5ex]{XI}
   &\mbox{b}&(\omega_\pi=1)
   &L(\nu^{1/2}\pi,\nu^{-1/2}\sigma)
   & \raisebox{2.5ex}[-2.5ex]{$\sigma^2$} 
   &\\\hline
 \end{array}
$$
\end{table}

\vspace{10ex}
\
%%%%%%%%

\begin{table}
\centering
\caption{The adjoint $L$-function $L(s,\Pi,\Ad)$}
\label{polestable}
\vspace{-3ex}
$$
\renewcommand{\arraystretch}{1.4}
 \begin{array}{|c|c|c|c|}
  \hline&&L(s,\Pi,{\rm Ad})&\mbox{ord}_{s=1}\\\hline\hline
% I
  &&L(s,1_{F^\times})^2L(s,\chi_1)L(s,\chi_1^{-1})L(s,\chi_2)L(s,\chi_2^{-1})&\\
  \raisebox{2.2ex}[-2.2ex]{I}&&L(s,\chi_1\chi_2)L(s,\chi_1^{-1}\chi_2^{-1})
   L(s,\chi_1\chi_2^{-1})L(s,\chi_1^{-1}\chi_2)&\raisebox{2.2ex}[-2.2ex]{$0$}\\\hline
% II
  &\mbox{a}&L(s,1_{F^\times})L(s,\chi^2)
   L(s,\chi^{-2})L(s,\nu)L(s,\chi^{-1}\nu^{1/2})L(s,\chi\nu^{1/2})&0\\\cline{2-4}
  {\rm II}&&L(s,1_{F^\times})^2L(s,\chi^2)L(s,\chi^{-2})L(s,\nu)L(s,\nu^{-1})&\\
   &\raisebox{2.2ex}[-2.2ex]{b}&L(s,\chi\nu^{-1/2})L(s,\chi^{-1}\nu^{1/2})
   L(s,\chi\nu^{1/2})L(s,\chi^{-1}\nu^{-1/2})&\raisebox{2.2ex}[-2.2ex]{$1$}\\\hline
% III
  &\mbox{a}&L(s,1_{F^\times})L(s,\nu)
   L(s,\nu\chi)L(s,\nu\chi^{-1})&0\\\cline{2-4}
  {\rm III}&&L(s,1_{F^\times})^2L(s,\chi)L(s,\chi^{-1})L(s,\nu)L(s,\nu^{-1})&\\
   &\raisebox{2.2ex}[-2.2ex]{b}
   &L(s,\chi\nu)L(s,\chi\nu^{-1})L(s,\chi^{-1}\nu)L(s,\chi^{-1}\nu^{-1})
   &\raisebox{2.2ex}[-2.2ex]{$1$ or $2$}\\\hline
% IV
  &\mbox{a}&L(s,\nu)L(s,\nu^3)&0\\\cline{2-4}
  &\mbox{b}&L(s,1_{F^\times})L(s,\nu)L(s,\nu^{-1})L(s,\nu^3)&1\\\cline{2-4}
  \raisebox{2.2ex}[-2.2ex]{IV}&
   \mbox{c}&L(s,1_{F^\times})L(s,\nu)L(s,\nu^{-1})L(s,\nu^2)
   L(s,\nu^3)L(s,\nu^{-3})&1\\\cline{2-4}
  &\mbox{d}&L(s,1_{F^\times})^2L(s,\nu)^2L(s,\nu^{-1})^2L(s,\nu^2)L(s,\nu^{-2})
   L(s,\nu^3)L(s,\nu^{-3})&2\\\hline
% V
  &\mbox{a}&L(s,\nu)^2L(s,\xi)L(s,\nu\xi)&0\\\cline{2-4}
  &\mbox{b}&L(s,1_{F^\times})L(s,\nu)^2L(s,\nu^{-1})
   L(s,\xi)L(s,\nu\xi)&1\\\cline{2-4}
  \raisebox{2.2ex}[-2.2ex]{V}&\mbox{c}&L(s,1_{F^\times})L(s,\nu)^2L(s,\nu^{-1})
   L(s,\xi)L(s,\nu\xi)&1\\\cline{2-4}
  &\mbox{d}&L(s,1_{F^\times})^2L(s,\nu)^2L(s,\nu^{-1})^2
   L(s,\xi)^2L(s,\nu\xi)L(s,\nu^{-1}\xi)&2\\\hline
% VI
  &\mbox{a}&&\\\cline{2-2}
  &\mbox{b}&\raisebox{2.2ex}[-2.2ex]{$L(s,1_{F^\times})L(s,\nu)^3$}
  &\raisebox{2.2ex}[-2.2ex]{0}\\\cline{2-4}
  \raisebox{2.2ex}[-2.2ex]{VI}&\mbox{c}
   &L(s,1_{F^\times})^2L(s,\nu)^3L(s,\nu^{-1})&1\\\cline{2-4}
  &\mbox{d}&L(s,1_{F^\times})^4L(s,\nu)^3L(s,\nu^{-1})^3&3\\\hline\hline
% VII
  {\rm VII}&&L(s,1_{F^\times})L(s,\pi,{\rm Ad}_{\GL(2)})
  L(s,\pi,{\rm Ad}_{\GL(2)}\otimes\chi)L(s,\pi,{\rm Ad}_{\GL(2)}\otimes\chi^{-1})
   &0\\\hline
% VIII
  &\mbox{a}&
   &\\\cline{2-2}
  \raisebox{2.2ex}[-2.2ex]{VIII}
  &\mbox{b}&\raisebox{2.2ex}[-2.2ex]{$L(s,1_{F^\times})
   L(s,\pi,{\rm Ad}_{\GL(2)})^3$}&\raisebox{2.2ex}[-2.2ex]{$0$}\\\hline
% IX
  &\mbox{a}&L(s,\xi)L(s,\pi,{\rm Ad}_{\GL(2)}\otimes\xi\nu)
  &0\\\cline{2-4}
  \raisebox{2.2ex}[-2.2ex]{IX}
  &\mbox{b}&L(s,1_{F^\times})L(s,\pi,{\rm Ad}_{\GL(2)})
  L(s,\pi,{\rm Ad}_{\GL(2)}\otimes\xi\nu)L(s,\pi,{\rm Ad}_{\GL(2)}\otimes\xi\nu^{-1})  
  &1\\\hline\hline
% X
  {\rm X}&&L(s,1_{F^\times})L(s,\pi,{\rm Ad}_{\GL(2)})L(s,\omega_\pi)L(s,\omega_\pi^{-1})&
   0\\\hline
% XI
  &\mbox{a}&L(s,\pi,{\rm Ad}_{\GL(2)})L(s,\nu)&0\\\cline{2-4}
  \raisebox{2.2ex}[-2.2ex]{XI}
   &\mbox{b}&L(s,1_{F^\times})L(s,\pi,{\rm Ad}_{\GL(2)})L(s,\nu)L(s,\nu^{-1})&1\\\hline
 \end{array}
$$
\end{table}

%\end{appendix}%
%
%

\end{document}